\documentclass[letterpaper, 10pt, conference]{ieeeconf}

\IEEEoverridecommandlockouts
\usepackage{xspace,amssymb,epsfig,syntonly}
\usepackage{epsfig,amsmath,color}
\usepackage{xspace,syntonly,empheq}
\usepackage{url}
\usepackage{bbm}
\usepackage{mathtools}
\usepackage{cite}
\usepackage{graphicx}
\usepackage{algorithm,algorithmic}
\usepackage{verbatim}
\usepackage{amssymb}
\usepackage{circuitikz}
\usepackage{enumerate}

\newcommand{\utwi}[1]{\mbox{\boldmath $#1$}}

\newcommand{\cD}{{\cal D}}

\newcommand{\cL}{{\cal{L}}}
\newcommand{\cN}{{\cal N}}

\newcommand{\cH}{{\cal H}}

\newcommand{\cX}{{\cal X}}

\newcommand{\ba}{{\bf a}}

\newcommand{\bh}{{\bf h}}

\newcommand{\br}{{\bf r}}

\newcommand{\bx}{{\bf x}}

\newcommand{\bw}{{\bf w}}

\newcommand{\bz}{{\bf z}}
\newcommand{\by}{{\bf y}}

\newcommand{\bT}{{\bf T}}

\newcommand{\blambda}{{\utwi{\lambda}}}

\newcommand{\bphi}{{\utwi{\phi}}}

\newcommand{\reals}{\mathbb{R}}

\newcommand{\sfI}{\textsf{I}}
\newcommand{\sfT}{\textsf{T}}
\newcommand{\sfF}{\textsf{F}}

\usepackage{epstopdf}

\newtheorem{proposition}{Proposition}

\newtheorem{theorem}{Theorem}
\newtheorem{corollary}{Corollary}

\newtheorem{assumption}{Assumption}
\newtheorem{remark}{Remark}

\title{\LARGE \bf On the Convergence of the Inexact Running \\ Krasnosel'ski\u{\i}-Mann Method}

\author{Emiliano Dall'Anese$^{1}$, Andrea Simonetto$^{2}$, Andrey Bernstein$^{3}$ 
\thanks{$^{1}$E. Dall'Anese is the University of Colorado Boulder; email: emiliano.dallanese@colorado.edu. $^{2}$A. Simonetto is with IBM Research Ireland; email: andrea.simonetto@ibm.com. $^{3}$A. Bernstein is with the National Renewable Energy Laboratory (NREL); email: andrey.bernstein@nrel.gov. The work of E. Dall'Anese was supported by NREL via APUP UGA-0-41026-109. Funds for A. Bernstein were provided by ARPA-e NODES. }
}

\begin{document}
  \maketitle

\begin{abstract}

This paper leverages a framework based on averaged operators to tackle the problem of tracking fixed points associated with maps that evolve over time. In particular, the paper considers the Krasnosel'ski\u{\i}-Mann  method in  a settings where: (i) the underlying map may change at each step of the algorithm, thus leading to a ``running'' implementation of the Krasnosel'ski\u{\i}-Mann method; and, (ii) an imperfect information of the map may be available. An imperfect
knowledge of 
the maps can capture cases where processors feature a finite precision or quantization errors, or the case where (part of) the map is obtained from measurements. The analytical results are applicable to inexact running algorithms for solving optimization problems, whenever the algorithmic steps can be written in the form of (a composition of) averaged operators; examples are provided for inexact running gradient methods and the forward-backward splitting method.  Convergence of the average fixed-point residual is investigated for the non-expansive case; 
linear convergence to a unique fixed-point trajectory is showed in the case of inexact running algorithms emerging from contractive operators.

\end{abstract}

\section{Introduction and Problem Formulation}
\label{sec:introduction}

The Banach-Picard method and its Krasnosel'ski\u{\i}-Mann (KM) variant have been leveraged to establish convergence of a number of iterative algorithmic frameworks for solving convex optimization problems as well as problems associated with (non)linear systems~\cite{bauschke2011convex,Combettes2002,monotonePrimer,Simonetto17,Mou15}. Focusing on the KM method, recall that an operator $\sfT: \cD \rightarrow \cD$, where $\cD$ is a nonempty convex subset of a finite-dimensional Hilbert space $\cH$ with a given norm $\|\cdot\|$, is non-expansive if it is $1$-Lipschitz in $\cH$; that is, $\forall \,  \bx, \by \in \cD$ one has that $\|\sfT(\bx) - \sfT(\by)\| \leq \|\bx - \by\|$. The KM algorithm involves the sequential application of the following operator starting from a point in $\cD$ (with $k$ the iteration index):
\begin{align}
\sfF_k := \sfI + \lambda_k
( \sfT - \sfI )  \label{eq:averaged_def} 
\end{align}
with $\sfI: \cH \rightarrow \cH$ the identity operator and $\{\lambda_k\}_{k \in \mathbb{N}}$ a sequence in $[0,1]$ satisfying $\sum_{k = 1}^{\infty} \lambda_k (1 - \lambda_k) = \infty$~\cite{bauschke2011convex}. Based on~\eqref{eq:averaged_def}, convergence of iterative algorithms for solving optimization problems can be cast as the problem of finding fixed points of a properly constructed non-expansive map $\sfT$ (which are also fixed points of $\sfF$). As another example, the operator-based representation~\eqref{eq:averaged_def} can be utilized to investigate convergence of  discrete-time linear systems~\cite{Belgioioso18}. 

The KM method~\eqref{eq:averaged_def} is known to converge weakly to a fixed point of $\sfT$~\cite{bauschke2011convex,cominetti2014rate,Belgioioso18,Themelis2016}; that is, taking the case of a constant value of $\lambda_k = \lambda $ as an example, one has that the average fixed-point residual of the map $\sfT$ after $K$ iterations can be bounded as~\cite{bauschke2011convex,cominetti2014rate}:
\begin{align}
\label{eq:aggregate_residual}
\frac{1}{K} \sum_{k = 1}^K \|\bx_k - \sfT(\bx_k)\|^2 \leq \frac{\|\bx_1 - \bx^*\|^2}{K \lambda (1 - \lambda)}
\end{align}
with $\bx^*$ a fixed point. See also the inexact~\cite{Bravo2017} and stochastic~\cite{combettes2017stochastic} variants, as well as more results on convergence of algorithms involving averaged non-expansive operators~\cite{Bauschke2015}. 

While~\eqref{eq:aggregate_residual} pertains
to problems where the map $\sfT$ is ``fixed'' during the execution of the KM algorithm and it is known, this paper revisits the convergence of the KM method in case of \emph{time-varying} and possibly \emph{inexact} maps. This setting is motivated by recent efforts to address the design and analysis of running algorithms for time-varying optimization problems~\cite{Simonetto_Asil14,Rahili2015,Simonetto17,OnlineFeedback}, with particular emphasis on feedback-based online optimization~\cite{OnlineFeedback,chen2018bandit}; additional works along these lines are in the context of  online optimization (see the representative works~\cite{Koppel2015,Hall2015,Shahrampour2018} and references therein) and learning  in dynamic environments~\cite{Dixit2018,Csaji2008}. In a time-varying optimization setting, the underlying cost, constraints, and problem inputs may change at every step (or a few steps) of the algorithm; therefore, pertinent tasks in this case involve the derivation of results for the \emph{tracking}  of optimal solution trajectories. Updates of the algorithms may be implemented inexactly due to finite-precision~\cite{zhang2018augmented} or because measurement feedback is utilized in lieu of model-based gradient computations~\cite{OnlineFeedback}. 
Counterparts of~\eqref{eq:aggregate_residual} are of interest for inexact running algorithms for problems with time-varying cost functions that are (locally) convex but not strongly convex; in case of problems with a (locally) strongly convex costs,  contractive arguments can be leveraged.

To concretely outline the problem, consider discretizing the temporal index as $t h$, $t \in \mathbb{N}$ and with $h$ a given interval (that will coincide with the time required to evaluate a map). Taking the normed space $(\mathbb{R}^m, \|\cdot\|)$ for the rest of the paper, consider a convex and closed set $\cD \subseteq \reals^m$ and a sequence of non-expansive mappings $\sfF_t: \cD \rightarrow \cD$. In particular, assume that $\sfF_t$ is $\alpha_t$-averaged; that is, it is a convex combination 
\begin{align}
\label{eq:averaged}
\sfF_t = (1 - \alpha_t) \sfI + \alpha_t \sfT_t
\end{align}
$\alpha_t \in (0,1)$. Starting from $\bx_1 \in \cD$, the \emph{running KM method} amounts to the execution of the following step at each $t$:  
\begin{align}
\bx_{t} = \sfF_t(\bx_{t-1}) = (1 - \alpha_t) \bx_{t-1} + \alpha_t \sfT_{t}(\bx_{t-1})  \,. \label{eq:running_KM} 
\end{align}
Different from the ``batch'' KM method -- especially when a Mann sequence $\{\lambda_k\}_{k \in \mathbb{N}}$ is utilized -- where~\eqref{eq:averaged_def} is executed within an interval $h$ until convergence,  the running algorithm~\eqref{eq:running_KM} boils down to a sequential application of time-varying $\alpha_t$-averaged maps. Preliminary results for the convergence of~\eqref{eq:running_KM} were provided in~\cite{Simonetto17}. 

The paper investigates the ability of the running algorithm~\eqref{eq:running_KM} to \emph{track} fixed points of the sequence of mappings $\{\sfT_{t}\}_{k \in \mathbb{N}}$, when in imperfect mapping $\hat{\sfT}_t:  \cD \rightarrow \cD$ is available. Notice that fixed points would be identified at each time $t$ only if the KM method~\eqref{eq:averaged_def} is executed to convergence at each $t$ (i.e., in a batch setting, instead of performing only one iteration) and the map $\sfT_{t}$ is known. This paper derives results similar to~\eqref{eq:aggregate_residual} for the inexact running KM method; results are also provided for the case of vanishing errors and vanishing fixed-point dynamics. The paper further considers  the case where the overall mappings $\{\sfF_t\}_{k \in \mathbb{N}}$ are contractions, and establishes linear convergence to the unique fixed-point trajectory. The proposed framework is then exemplified for inexact running projected gradient  and forward-backward splitting methods for solving time-varying convex optimization problems. Overall, the paper provides contributions over our previous work~\cite{Bernstein_fixedpoint} on running Banach-Picard method, where linear convergence results where established in case of time-varying contractive maps, possibly corrupted by errors.
Stochastic time-varying-fixed problems were considered in~\cite[Th.~20]{Csaji2008}; here, we focus on bounded errors on averaged operators, and leave stochastic errors as a follow on research opportunity.

\section{Inexact Running Algorithm}
\label{sec:inexact}

Let $\bx^{\star}_t$ be a fixed point of the self-mapping $\sfF_t$; that is, $\bx^{\star}_t = \sfF_t(\bx^{\star}_t)$. 
If the vectors $\{\bx^{\star}_t\}_{t \in \mathbb{N}}$ satisfy the equation $\bx^{\star}_t = \sfF_t(\bx^{\star}_t)$ for each $t \in \mathbb{N}$, then we refer to $\{\bx^{\star}_t\}_{t \in \mathbb{N}}$ as a sequence of fixed points. If the mappings $\{\sfF_t\}_{t \in \mathbb{N}}$ are averaged, multiple sequences  $\{\bx^{\star}_t\}_{t \in \mathbb{N}}$ may exist; since $\sfF_t = (1 - \alpha_t) \sfI + \alpha_t \sfT_{t}$, $\bx^{\star}_t$ is also a fixed point of $\sfT_t$.
When $\{\sfF_t\}_{t \in \mathbb{N}}$ are contractions, only one sequence exists by the Banach fixed-point theorem.  To characterize the variability of a fixed-point sequence, we assume that there exists \emph{a} sequence of fixed points $\{\bx^{\star}_t\}_{t \in \mathbb{N}}$, for which there exists a finite and non-negative  sequence of scalars $\{\sigma_t\}_{t \in \mathbb{N}}$, such that
\begin{align}
 \|\bx^{\star}_{t+1} - \bx^{\star}_{t}\| \leq \sigma_t \, . \label{eq:sigma} 
\end{align}
for all $t$. If $\sfF_{t+1} = \sfF_t$ then one has that $\sigma_t = 0$, and we are recover the time-invariant case.

Consider now a mapping $\hat{\sfT}_t: \cD \rightarrow \cD$, which is an approximation of $\sfT_t$ in the following sense.

\vspace{.1cm}

\begin{assumption}[Bounded approximation error] \label{asm:assump_error}
For each $t \in \mathbb{N}$ and for all $\bx \in \cD$, it holds that $\hat{\sfT}_t(\bx) \in \cD$. Further, there exists a scalar $e_{\sfT,t} < + \infty$ such that 
\begin{align}
\max_{\bx \in \cD}\|\sfT_t(\bx) - \hat{\sfT}_t(\bx)   \| \leq e_{\sfT,t} \, .
 \label{eq:bound_errro_g} 
\end{align}
\end{assumption}

\vspace{.1cm}

The condition~\eqref{eq:bound_errro_g} simply asserts that the error in the map is bounded; it can be deterministic or stochastic (and i.i.d over time), but with finite support. Accordingly, define the approximate $\alpha_t$-averaged map $\hat{\sfF}_t$ as: 
\begin{align}
\hat{\sfF}_t(\bx) := (1 - \alpha_t) \bx + \alpha_t \hat{\sfT}_t(\bx)  \label{eq:averaged_approximate} 
\end{align}
Based on~\eqref{eq:averaged_approximate}, and given an initial point $\bx_1 \in \cD$ the \emph{inexact running KM} algorithm is given by [cf.~\eqref{eq:running_KM}]:    
\begin{align}
\bx_{t} = \hat{\sfF}_t(\bx_{t-1}) = (1 - \alpha_t) \bx_{t-1} + \alpha_t \hat{\sfT}_t(\bx_{t-1}) \, . 
 \label{eq:running_mk_err} 
\end{align}
In the next section, tracking of a sequence of fixed points $\{\bx^{\star}_t\}_{t \in \mathbb{N}}$ via~\eqref{eq:running_mk_err} will be investigated.

\section{Convergence}
\label{sec:convergence}
 
This section will characterize the performance of the inexact running KM method in two different settings: 

\noindent \emph{i)} The map $\sfT_t$ is non-expansive and $\sfF_t$ is $\alpha_t$-averaged; and, 

\noindent \emph{ii)} The map $\sfF_t$ is a contraction. 

\noindent   It is worth pointing out that for generic non-expansive maps, the sequence generated by the Banach-Picard iteration may fail to produce a fixed point even in a static case; the structure of~\eqref{eq:running_mk_err} will however facilitate the derivation of convergence results. Regarding the second case, notice that if $\sfT_t$ is contractive then $\sfF_t$ is contractive; however, the converse is not necessarily true. We start by outlining the following standard assumptions~\cite{bauschke2011convex,cominetti2014rate}.

 
 

\vspace{.1cm}

\begin{assumption}[Lipshitz maps] \label{asm:assump_L}
There exists a scalar $0 \leq L_t \leq 1$ such that $\| \sfF_t(\bx) - \sfF_t(\bx') \| \leq L_t \|\bx - \bx' \|$  for all $\bx, \bx' \in \cD$.
\end{assumption}

\vspace{.1cm}

\begin{assumption}[Bounded maps] 
\label{asm:assump_bounded_f}
There exists a scalar $M_t < + \infty$ such that 
\begin{align}
\max_{\bx \in \cD}\left\|\sfF_t(\bx) \right\| \leq M_{t} \, , \,\,\,\, \max_{\bx \in \cD}\|\hat{\sfF}_t(\bx) \| \leq M_{t} .
\end{align}
\end{assumption}

\vspace{.1cm}

\noindent If $\cD$ is compact, then $M_{t}$ can be taken, in the worst case, to be the radius of $\cD$. For subsequent developments, define $M := \sup_{t} \{M_t\}$, $\sigma := \sup_{t} \{\sigma_t\}$, $e_{\sfT} := \sup_{t} \{e_{\sfT,t}\}$, and $\alpha := \sup_{t} \{\alpha_t\}$. The following result pertains to the case where  $\sfF_t$ is $\alpha_t$-averaged.

\vspace{.1cm}

\begin{theorem}
\label{thm:residual}
Consider a sequence of $\alpha_t$-averaged operators $\sfF_t = (1 - \alpha_t) I  + \alpha_t \sfT_t$, $t = 1, \ldots, T$, and assume that there exists a sequence of vectors $\{\bx^{\star}_t\}_{t = 1}^T$ that satisfy the equation $\bx^{\star}_t = \sfF_t(\bx^{\star}_t)$ for each $t = 1, \ldots, T$. Suppose that Assumptions~\ref{asm:assump_error}--\ref{asm:assump_bounded_f} hold, and take $\bx_1 \in \cD$. Then, the following bound holds for the algorithm~\eqref{eq:running_mk_err}:
\begin{align}
&\sum_{t = 1}^T \alpha_t(1 - \alpha_t) \|\bx_t - \sfT_t(\bx_t)\|^2 \leq \|\bx_1 - \bx_1^{\star}\|^2 +  \sum_{t = 1}^T r_t
 \label{eq:fp_residual} 
\end{align}
where $r_t :=  \alpha_t e_{\sfT,t} (4 M_t   +  \alpha_t e_{\sfT,t}) +  \sigma_t (4 M_t + \sigma_t)$. In particular, one has that:
\begin{align}
& \hspace{-.1cm} \frac{1}{T}\sum_{t = 1}^T \alpha_t(1 - \alpha_t) \|\bx_t - \sfT_t(\bx_t)\|^2 \leq \frac{1}{T} \|\bx_1 - \bx^{\star}_1\|^2 \hspace{-.1cm} +  r
 \label{eq:fp_residual_2} 
\end{align}
\begin{align}
& \hspace{-.5cm} \frac{1}{T}\sum_{t = 1}^T \frac{1 - \alpha_t}{\alpha_t} \|\bx_t - \sfF_t(\bx_t)\|^2 \leq \frac{1}{T} \|\bx_1 - \bx^{\star}_1\|^2 +  r
 \label{eq:fp_residual_3} 
\end{align}
with $r := \alpha e_{\sfT} (4 M   +  \alpha e_{\sfT}) +  \sigma (4 M + \sigma) $. 
\end{theorem}

\vspace{.1cm}

\emph{Proof.} See Appendix~\ref{sec:proof_thm_residual}

\vspace{.1cm}

Bounds~\eqref{eq:fp_residual_2}--\eqref{eq:fp_residual_3} imply convergence in mean of the fixed-point residual to a ball centered at $0$; the size of the ball
depends on the bound on the variability of the fixed-point trajectories, on the size of the image of the operators, and on the approximation errors for the maps. An immediate follow-up from~\eqref{eq:fp_residual_2}--\eqref{eq:fp_residual_3} is the following asymptotic result:
\begin{align}
& \limsup_{T\rightarrow + \infty}  \frac{1}{T}\sum_{t = 1}^T  \|\bx_t - \sfF_t(\bx_t)\|^2 \leq  r \bar{\alpha}^{-1}
 \label{eq:fp_residual_asymp} 
\end{align}
where $\bar{\alpha} := \inf_{t = 1, \ldots T} \{(1 - \alpha_t)/\alpha_t\}$. A similar result can be derived for the mean  of $\|\bx_t - \sfT_t(\bx_t)\|^2$.

It is worth pointing out that, when $\sigma = 0$, the bound in~\eqref{eq:fp_residual_asymp} reduces to $\alpha e_{\sfT} (4 M   +  \alpha e_{\sfT}) \bar{\alpha}^{-1}$, and the bounds therefore capture the effect of the approximate maps. In case of perfect mappings,~\eqref{eq:fp_residual_asymp} boils down to~\eqref{eq:aggregate_residual}~\cite{bauschke2011convex,cominetti2014rate}. Motivated by this, the next results will deal with vanishing errors and fixed-point dynamics, which is increasingly motivated by learning in bandit settings (where the maps are learned online while the algorithm is running).  

\vspace{.1cm} 

\begin{corollary}
\label{cor:vanishing_error}
Suppose\footnote{A relation $f(n) = o(g(n))$ signifies that for every positive constant $\varphi$ there exists $N$ such that $|f(n)|\leq \varphi |g(n)|$ for all $n\geq N$.} that for each $T$, one has that
\begin{align}
\label{eq:sum_errors}
\sum_{t = 1}^T  e_{\sfT,t} = o(T) , 
\end{align}
i.e., $\sum_{t = 1}^T  e_{\sfT,t}$ grows sublinearly in $T$. If Assumptions~\ref{asm:assump_error}--\ref{asm:assump_bounded_f} hold,  then, for the algorithm~\eqref{eq:running_mk_err}, the fixed-point residual $\|\bx_t - \sfT_t(\bx_t)\|$ 
 converges to:
\begin{align}
& \limsup_{T \rightarrow \infty}\frac{1}{T}\sum_{t = 1}^T  \|\bx_t - \sfT_t(\bx_t)\|^2 \leq \check{\alpha}^{-1} \sigma (4 M + \sigma) 
 \label{eq:fp_residual_vani}  \\
& \limsup_{T \rightarrow \infty}\frac{1}{T}\sum_{t = 1}^T  \|\bx_t - \sfF_t(\bx_t)\|^2 \leq \bar{\alpha}^{-1} \sigma (4 M + \sigma) 
 \label{eq:fp_residual_vani2}  
\end{align}
where $\check{\alpha} := \inf_{t = 1, \ldots T} \{\alpha_t(1 - \alpha_t)\}$. 
\end{corollary}

\vspace{.1cm} 

\emph{Proof.} See Appendix~\ref{sec:proof_cor_vanishing_error}.  

\vspace{.1cm} 

\begin{corollary}
\label{cor:vanishing_error}
Suppose that for each $T$, one has that
\begin{align}
\label{eq:sum_dynamics}
\sum_{t = 1}^T  \sigma_t = o(T) 
\end{align}
i.e., $\sum_{t = 1}^T \sigma_t$ grows sublinearly in $T$. Assume further that~\eqref{eq:sum_errors} holds. Then,  under  Assumptions~\ref{asm:assump_error}--\ref{asm:assump_bounded_f}, for the algorithm~\eqref{eq:running_mk_err} one has that $\lim_{t \rightarrow \infty} \|\bx_t - \sfT_t(\bx_t)\|^2 = 0$ and $\lim_{t \rightarrow \infty} \|\bx_t - \sfF_t(\bx_t)\|^2 = 0$. 
\end{corollary}

\vspace{.1cm}

For completeness, we now turn the attention to  convergence results for contractive operators. The following holds. 

\vspace{.1cm}

\begin{theorem}
\label{thm:convergence}
Consider a sequence of contractive mappings of the form $\sfF_t = (1 - \alpha_t) \sfI  + \alpha_t \sfT_t$, $t = 1, \ldots, T$ and let $\{\bx^{\star}_t\}_{t = 1}^T$ be the trajectory of  fixed points. Let $\{\bx_t\}_{t = 1}^T$ be a sequence generated by the algorithm~\eqref{eq:running_mk_err}, with  $\bx_1 \in \cD$. Suppose that Assumptions~\ref{asm:assump_error}--\ref{asm:assump_bounded_f} hold. 
Then, at each time $t$, it holds that:
\begin{align}
\|\bx_{t+1} - \bx^{\star}_{t+1}\| & \leq c^{(t,0)} \|\bx_1 - \bx^{\star}_1\| \nonumber \\
& + \sum_{\tau = 1}^t c^{(t, \tau)} \left(\alpha_\tau e_{\sfT,\tau} + \sigma_\tau  \right)  \label{eq:bound_iter} 
\end{align}
for each $t$, where
\begin{equation} \label{eqn:beta}
c^{(t, \tau)} := 
\begin{cases}
\prod_{\ell = \tau + 1}^t L_{\ell}, & \textrm{if~} \tau = 0, \ldots, t - 1 \\
1, &  \textrm{if~} \tau = t.
\end{cases}
\end{equation}
Suppose further that Assumption~\ref{asm:assump_L} holds with $L_t < 1$ for all $t$. Then, $\{\bx^{\star}_t\}_{t = 1}^T$ is unique and the following asymptotic bound holds for the algorithm~\eqref{eq:running_mk_err}:
\begin{align}
& \lim_{t \rightarrow  \infty} \sup \|\bx_t - \bx^{\star}_t\| \leq 
\frac{\gamma}{1-L} 
 \label{eq:fp_linear} 
\end{align}
where $\gamma := \alpha e_{\sfT} + \sigma$ and $L := \sup_t \{L_t\}$.  
\end{theorem}

\vspace{.1cm}

\emph{Proof}. See Appendix~\ref{sec:proof_thm_convergence}.

\vspace{.1cm}

Bound~\eqref{eq:fp_linear} in similar to~\cite{Bernstein_fixedpoint}, but customized for the operators considered here. In case of vanishing errors and dynamics, the following results readily hold.

\begin{corollary}
Suppose that~\eqref{eq:sum_errors} holds. Then, if Assumption~\ref{asm:assump_L} holds with $L_t < 1$ for all $t$, then 
\begin{align}
& \lim_{t \rightarrow  \infty} \sup \|\bx_t - \bx^{\star}_t\| \leq 
\frac{\sigma}{1-L} \, .
 \label{eq:fp_linear_2} 
\end{align} 
Additionally, if~\eqref{eq:sum_dynamics} holds, then $\lim_{t \rightarrow \infty} \sup \|\bx_t - \bx^{\star}_t\| = 0$. 
\end{corollary}

\begin{remark}
When a predictable sequence is available, one could reduce the error ball $r$ to a sublinear function of $T$ by properly tuning the sequence $\{\alpha_t\}$, even if $\sigma_t$ and $e_{\sfT}$ do not vanish; see, for example, the framework in~\cite{Jadbabaie2015} for adaptive optimistic mirror descent methods.
Due to space limitations, we leave the derivation of these results for future efforts.
\end{remark}

\vspace{.1cm}

\begin{remark}
Proof techniques in~\cite{Berinde2007} presuppose particular sequences $\{\alpha_t\}$ and $\{e_{\sfT,t}\}$ to  establish convergence results for e.g., static non-expansive and strictly pseudocontractive  maps  (see, e.g.,  Theorems 6.1 and 6.2) as well as for (static) maps defined in Banach spaces (see, e.g., Theorem 6.8). Adopting the sequences $\{\alpha_t\}$ in~\cite{Berinde2007} might not be possible in a time-varying setting, especially when $\alpha_t \rightarrow 0$ for $t \rightarrow \infty$; however, future efforts will look at possible extensions of the techniques in~\cite{Berinde2007}  in the time-varying case.
 \end{remark}

\section{Examples of applications}
\label{sec:examples}

The objective of this  section is to show that a number of inexact running algorithms for time-varying optimization problems can be analyzed by leveraging the operator-based framework proposed in this paper. In particular, this section focuses on inexact running gradient methods and forward-backward splitting algorithms. Additional applications are possible~\cite{monotonePrimer}, but are not included due to space limitations.   

\subsection{Running gradient method with errors}
\label{sec:gradient}

Recall that the temporal index is discretized as $t h$, $t \in \mathbb{N}$,  with $h$ a given interval (that can coincide with the time required to perform one algorithmic step). Consider the following time-varying optimization problem

\begin{align}
\label{eq:gradient_1}
 (\textrm{P1}_{t}): \quad \min_{\bx \in \cX_t} f_t(\bx)
\end{align}
where $f_t: \mathbb{R}^n  \rightarrow \mathbb{R}$ is a convex, closed, and proper (CCP) function at each time $t$, and $\cX_t$ is a convex and compact set at each time $t$. Assume that $f$ is strongly smooth with parameter $K_t > 0$.  
Notice that solving the problem~\eqref{eq:gradient_1} is equivalent to finding the zeros of $\nabla f_t + \cN_{\cX_t}$, where $\cN_{\cX_t}$ is the  normal cone operator for the set $\cX_t$.

A running version of the projected gradient method for solving~\eqref{eq:gradient_1} is given by:
\begin{align}
\bx_{t+1} = \textrm{proj}_{\cX_t}\{\bx_{t} - \nu \nabla f_t(\bx_t)\}
\label{eq:gradient_2}
\end{align}
for a given step size $\nu > 0$. Let  $\by_t$ be a measurement or an estimate of the gradient $\nabla f_t(\bx_t)$; then, an inexact running projected gradient method is given by:
\begin{align}
\bx_{t+1} = \textrm{proj}_{\cX_t}\{\bx_{t} - \nu \by_t\}
\label{eq:gradient_inexact} \, .
\end{align}
In this setting, the bounds~\eqref{eq:fp_residual} and~\eqref{eq:fp_residual_2} will be utilized to derive tracking results for~\eqref{eq:gradient_inexact} for the case where the function $f_t$ is convex, but not strongly convex; on the other hand,~\eqref{eq:fp_linear} will be utilized for the case where $f_t$ is strongly convex uniformly in time.

For simplicity, focus first on the case where $\cX_t = \mathbb{R}^m$. Take $\nu \in (0, 2/K)$, with $K := \sup_t \{K_t\}$, so that the operator $\sfI - \nu \nabla f_t$ is averaged; that is,
\begin{align}
\sfI - \nu \nabla f_t = \left(1 - \nu K_t/2 \right) \sfI + {\nu K_t}/{2} \left(\sfI - 2/K_t \nabla f_t  \right)
\label{eq:gradient_3} 
\end{align}
which is in the form of~\eqref{eq:averaged} with $\alpha_t = \nu K_t/2$ and $\sfT_t = \sfI - \frac{2}{K_t} \nabla f_t$~\cite{Combettes15}. On the other hand, the approximate map $\hat{\sfT}_t$ is given by $\hat{\sfT}_t(\bx_t) = \bx_t - \frac{2}{K_t} \by_t$.  
Therefore, for the case where $\cX_t = \mathbb{R}^m$, one has that:
\begin{align}
\|\hat{\sfT}_t(\bx) - \sfT_t(\bx)\| & \leq 2 K_t^{-1} \|\nabla f_t(\bx) - \by\| \, .
\label{eq:gradient_error} 
\end{align}
Therefore, if there exists scalar $e_{y,t} < + \infty$ so that $ \|\nabla f_t(\bx) - \by\| \leq e_{y,t}$~\cite{OnlineFeedback}, $e_{\sfT,t}$ in~\eqref{eq:bound_errro_g} amounts to: 
\begin{align}
e_{\sfT,t} = 2 K_t^{-1} e_{y,t} .
\label{eq:gradient_error_2} 
\end{align}

The results for the inexact running projected gradient method are presented in the following proposition.   

\vspace{.1cm} 

\begin{proposition}
\label{prop:gradient}
Let $\nu \in (0, 2/K)$, and let $\{\bx_{t}\}$ be a sequence generated by~\eqref{eq:gradient_inexact}. Assume that there exists scalar $e_{y,t} < + \infty$ so that $ \|\nabla f_t(\bx) - \by\| \leq e_{y,t}$. Then, one has that~\eqref{eq:gradient_inexact} is an inexact averaged operator with $\alpha_t = 1/(2 - \nu K_t/2)$ and 
\begin{align}
e_{\sfT,t} =  \left(2 \nu - {\nu^2 K_t}/{2} \right) e_{y,t} \, .
\label{eq:gradient_error_3} 
\end{align}
For the algorithm~\eqref{eq:gradient_2}:

\noindent (i) The bounds~\eqref{eq:fp_residual},~\eqref{eq:fp_residual_2}, and~\eqref{eq:bound_iter} hold with  $e_{\sfT,t}$ as in~\eqref{eq:gradient_error_3};

\noindent (ii) Suppose further that $f_t$ is strongly convex with constant $k_t$; then,~\eqref{eq:fp_linear} hold with $L_t = \min\{|1 - \nu k_t|,|1 - \nu K_t|\}$.    

\end{proposition}

\vspace{.1cm} 

\emph{Proof.} See Appendix~\ref{sec:proof_prop_gradient}. 


\subsection{Inexact forward-backward splitting method }
\label{sec:forward-backward}

Consider the following time-varying problem~\cite{Dixit2018}
\begin{align}
\label{eq:forward-backward_1} 
(\textrm{P2}_{t}):\quad \min_{\bx \in \cX_t} f_t(\bx) + g_t(\bx) 
\end{align}
where $f_t: \mathbb{R}^n  \rightarrow \mathbb{R}$ and $g_t: \mathbb{R}^n \rightarrow \mathbb{R}$ are CCP functions at each time $t$, and $\cX_t$ is a convex and compact set at each time $t$. Assume that $f_t$ is strongly smooth with parameter $K_t > 0$ for all $t$, and suppose that $g_t$ is not differentiable.

A running version of the forward-backward splitting method for solving~\eqref{eq:forward-backward_1} is given by:
\begin{align}
\bx_{t+1} = \textrm{prox}_{g_t, \cX_t, \nu}\{\bx_{t} - \nu \nabla f_t(\bx_t)\}
\label{eq:forward-backward_2}
\end{align}
where 
\begin{align}
\textrm{prox}_{g_t, \cX_t, \nu}\{\by\} := \arg \min_{\bx \in \cX_t} \left\{g_t(\bx) + \frac{1}{2 \nu} \|\bx - \by\|_2^2\right\}
\label{eq:forward-backward_3}
\end{align}
is the proximal operator. If $\nu \in (0, 2/K)$, then the update~\eqref{eq:forward-backward_2} is given by the composition of a proximal operator and the operator $\sfI - \nu \nabla f_t$. The proximal operator is $\frac{1}{2}$-averaged~\cite{Combettes15,monotonePrimer}, whereas $\sfI - \nu \nabla f_t$ is an averaged operator with $\alpha_t = \nu K_t/2$, whenever $\nu \in (0, 2/K)$. Therefore, since the composition of averaged operators is an averaged operator, if follows from~\cite{Combettes15} that~\eqref{eq:forward-backward_2} is an averaged operator with $\alpha_t = 1/(2 - \nu K_t/2)$. 

An inexact version of the running forward-backward splitting method for solving~\eqref{eq:forward-backward_1} is given by:
\begin{align}
\bx_{t+1} = \textrm{prox}_{g_t, \cX_t, \nu}\{\bx_{t} - \nu \by_t\}
\label{eq:forward-backward_4}
\end{align}
where $\by_t$ is a measurement or an estimate of $\nabla f_t(\bx_t)$. Assuming that there exists scalar $e_{y,t} < + \infty$ so that $ \|\nabla f_t(\bx) - \by\| \leq e_{y,t}$,  results similar to Proposition~\ref{prop:gradient} apply to the inexact running forward-backward splitting method~\eqref{eq:forward-backward_4}. In particular,~\eqref{eq:fp_residual} and~\eqref{eq:fp_residual_2} bound the tracking error for~\eqref{eq:forward-backward_4} when the function $f_t$ is not strongly convex.

\section{Illustrative Numerical Results}
\label{sec:results}

\begin{figure}[t!]
  \centering
  \includegraphics[width=.4\columnwidth]{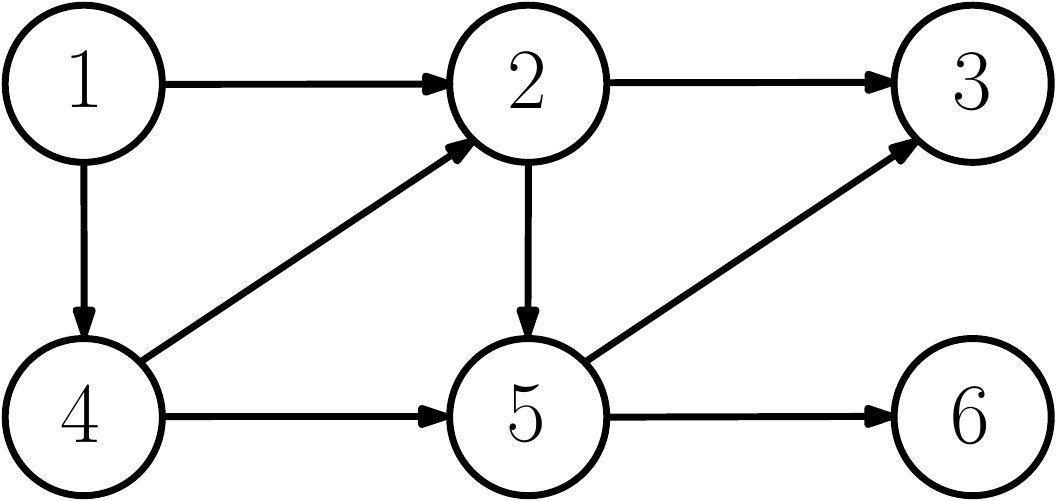}
\vspace{-.3cm}
\caption{Network utilized in the numerical results}
\label{fig:F_comm_net}
\vspace{-.2cm}
\end{figure}

As an illustrative example, we consider the network in Fig.~\ref{fig:F_comm_net} with 6 nodes and 8 links. The routing matrix is based on the directed edges. Let $z(i,s)$ denote the  rate generated at node $i$ for traffic $s$ and $r(ij,s)$ the flow between noted $i$ and $j$ for traffic $s$. consider then the following problem:
\label{eqn:sampledProblemComm}
\begin{align} 
&\min_{\substack{\{\bz, \br \}  \in \cX_t } }\sum_{i,s} - \kappa(i,s) \log(1+z(i,s)) + \ba_t^\sfT \br \label{eq:obj_p0_Comm} 
\end{align} 
where $\bz$ and $\br$ stack the traffic rates and link rates for brevity, $\kappa(i,s)$ and $\ba$ are given positive coefficients and the set $\cX_t$ is built based on: i) the flow-conservation constraints $\bz_s = \bT (\br^s + \bw^s_t)$ per flow $s$, where $\bT$ is the routing matrix and $\bw^s_t$ is a time-varying exogenous flow (of uncontrollable traffic); ii) the per-link capacity constraints, where the capacity of link $(i,j)$ is given by $\log(1+ p(i,j) h(i,j))$, with $p(i,j)$ the transmit power and $h(i,j)$ the normalized channel gain; and, iii) the non-negativity constraints on the traffic rates. Assume that two traffic flows are generated by nodes $1$ and $4$, and they are received at nodes $3$ and $6$, respectively. 

We utilize~\eqref{eq:gradient_inexact}. Errors and time variability of the problem are introduced as follows: 

\noindent $\bullet$ Gradient errors: the gradient of the cost $\kappa(i,s) \log(1+z(i,s))$ for each exogenous traffic flow is estimated using a multi-point  bandit feedback~\cite{AgarwalXiao10,chen2018bandit}; the estimation error depends on the number  of functional evaluations in constructing the proxy of the gradient in~\eqref{eq:gradient_inexact}. 

\begin{figure}[t!]
  \centering
  \includegraphics[width=1.0\columnwidth]{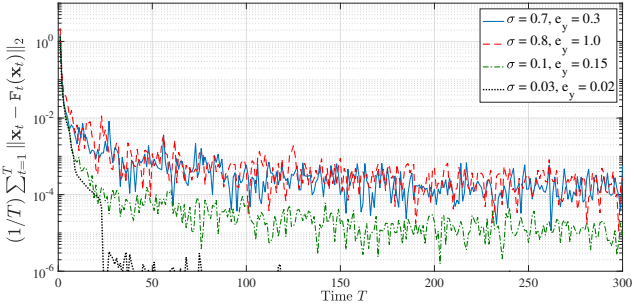}
\vspace{-.6cm}
\caption{Evolution of the fixed-point residual.}
\label{fig:F_fp_residual}
\vspace{-.6cm}
\end{figure}

\noindent $\bullet$  Solution dynamics: at each time step, the channel gain of links are generated by using a complex Gaussian random variable with mean $1 + \jmath 1$ and a given variance $v_c$ for both real and imaginary parts; the transmit power for each node is a  Gaussian random variable with mean $1$ and a variance $v_p$; the exogenous traffics are random with mean $[0.2, 0.3, 0.3, 0.4, 0.5, 0.2, 0.1, 0.4]$ and a given variance; and, the cost is perturbed by modifying $\ba_t$. Different values for $\sigma_t$ and $\sigma$ are obtained by varying the variance of these random variables. 
Figure~\ref{fig:F_fp_residual} illustrates the evolution of the fixed-point residual $(1/T)\sum_{t = 1}^T  \|\bx_t - \sfF_t(\bx_t)\|^2 $, for different values of $\sigma$ and the normalized error in the gradient estimate $e_y$. Optimal  rates are in the order of $0.6 - 1.7$; $\sigma = 0.7$ implies a $20\%$ worst-case variation in the solution between consecutive time steps, while $\sigma = 0.03$ leads to a $1\%$ variation. It can be seen that the fixed-point residual flattens, with an error that increase with the increasing of  $\sigma$ and  $e_y$, thus corroborating the proposed analytical results. 

\appendix

\subsection{Proof of Theorem~\ref{thm:residual}}
\label{sec:proof_thm_residual}

Consider $\|\bx_{t+1} - \bx^{\star}_{t+1}\|^2$, which can be bounded as follows by using the definition of $\sigma_t$:
\begin{subequations}
\begin{align}
\|\bx_{t+1} - \bx^{\star}_{t+1}\|^2 & = \|\bx_{t+1} - \bx^{\star}_{t} - (\bx^{\star}_{t+1} - \bx^{\star}_{t})\|^2 \\
& = \|\bx_{t+1} - \bx^{\star}_{t}\|^2 + \|\bx^{\star}_{t+1} - \bx^{\star}_{t}\|^2 \nonumber \\
& \hspace{.4cm} - 2 (\bx_{t+1} - \bx^{\star}_{t})^\sfT (\bx^{\star}_{t+1} - \bx^{\star}_{t}) \\
& \leq \|\bx_{t+1} - \bx^{\star}_{t}\|^2 + \sigma_t^2 \nonumber \\
& \hspace{.4cm} - 2 (\bx_{t+1} - \bx^{\star}_{t})^\sfT (\bx^{\star}_{t+1} - \bx^{\star}_{t}) \, . ~\label{eq:proof_afp_1}
\end{align}
\end{subequations}

The term $\|\bx_{t+1} - \bx^{\star}_{t}\|^2$ can be expanded as:
\begin{subequations}
\begin{align}
 \|\bx_{t+1} - \bx^{\star}_{t}\|^2 & = \|\hat{\sfF}_t(\bx_{t})  - \bx^{\star}_{t}\|^2 \\
&  = \|(1-\alpha_t)\bx_{t} + \alpha_t \hat{\sfT}_t(\bx_{t})  - \bx^{\star}_{t}\|^2 \\
& = \|(1-\alpha_t)\bx_{t} + \alpha_t \sfT_t(\bx_{t})  - \bx^{\star}_{t} \nonumber \\
& \hspace{1.2cm} + \alpha_t(\hat{\sfT}_t(\bx_{t}) - \sfT_t(\bx_{t}))\|^2 \, .  \label{eq:proof_afp_2}
\end{align}
\end{subequations}
Let $\ba_t := (1-\alpha_t)\bx_{t} + \alpha_t \sfT_t(\bx_{t})  - \bx^{\star}_{t}$ for brevity. Then,~\eqref{eq:proof_afp_2} can be further bounded as:

\begin{subequations}
\begin{align}
& \|\ba_t + \alpha_t(\hat{\sfT}_t(\bx_{t}) - \sfT_t(\bx_{t})) \|^2 \nonumber \\
& \hspace{.7cm} \leq \|\ba_t\|^2 + \alpha_t^2\|\hat{\sfT}_t(\bx_{t}) - \sfT_t(\bx_{t})\|^2 \nonumber \\
& \hspace{1.5cm} + 2 \alpha_t \ba_t^\sfT (\hat{\sfT}_t(\bx_{t}) - \sfT_t(\bx_{t})) \label{eq:proof_afp_3} \\
& \hspace{.7cm} \leq \|\ba_t\|^2 + \alpha_t^2\|\hat{\sfT}_t(\bx_{t}) - \sfT_t(\bx_{t})\|^2 \nonumber \\
& \hspace{1.5cm} + 2 \alpha_t \|\ba_t\| \|\hat{\sfT}_t(\bx_{t}) - \sfT_t(\bx_{t})\|  \label{eq:proof_afp_3} \\
& \hspace{.7cm} \leq \|\ba_t\|^2 + \alpha_t^2 e_{\sfT,t}^2  + 2 \alpha_t \|\ba_t\| e_{\sfT,t} \label{eq:proof_afp_4} 
\end{align}
\end{subequations}
To bound $\|\ba_t\|^2$, consider the following inequality, valid for any vectors $\bx \in \mathbb{R}^2$, $\by \in \mathbb{R}^2$ and scalar $\theta$: 
\begin{align}
\|(1-\theta) \bx + \theta \by\|^2 & = (1-\theta) \|\bx\|^2 + \theta \|\by\|^2 \nonumber \\
& \hspace{.4cm} - \theta (1-\theta) \|\bx - \by\|^2 \, .
\label{eq:proof_afp_5} 
\end{align}
Then, using~\eqref{eq:proof_afp_5} and the fact that $\bx^{\star}_{t} = (1 - \alpha_t) \bx^{\star}_{t} + \alpha_t \sfT_t(\bx^{\star}_{t})$, one has that:
\begin{subequations}
\begin{align}
& \|\ba_t\|^2  = \|(1-\alpha_t)\bx_{t} + \alpha_t \sfT_t(\bx_{t})  - \bx^{\star}_{t}\|^2 \\
& = \|(1-\alpha_t)(\bx_{t} - \bx^{\star}_{t}) + \alpha_t(\sfT_t(\bx_{t}) - \sfT_t(\bx^{\star}_{t})) \|^2 \\
& = (1-\alpha_t)\|\bx_{t} - \bx^{\star}_{t}\|^2 + \alpha_t \|\sfT_t(\bx_{t}) - \sfT_t(\bx^{\star}_{t})\|^2 \nonumber \\
& \hspace{.4cm} - \alpha_t(1 - \alpha_t) \|\sfT_t(\bx_{t}) - \bx_{t}\|^2 \label{eq:proof_afp_6} \\
& \leq \|\bx_{t} - \bx^{\star}_{t}\|^2 - \alpha_t(1 - \alpha_t) \|\sfT_t(\bx_{t}) - \bx_{t}\|^2 \label{eq:proof_afp_7}
\end{align}
\end{subequations}
where the non-expansiveness of $\sfT_t$ was used to obtain~\eqref{eq:proof_afp_7}. To bound $\|\ba_t\|$, it follows from Assumption~\ref{asm:assump_bounded_f} that: 
\begin{subequations}
\begin{align}
 \|\ba_t\| & = \|(1-\alpha_t)\bx_{t} + \alpha_t \sfT_t(\bx_{t})  - \bx^{\star}_{t}\| \\
& = \|\sfF_t(\bx_{t-1}) - \sfF_t(\bx^{\star}_{t})\|
\leq 2 M_t  \label{eq:proof_afp_9}
\end{align}
\end{subequations}
Regarding the third term on the right-hand-side of~\eqref{eq:proof_afp_1}, one can show that:
\begin{subequations}
\begin{align}
& - (\bx_{t+1} - \bx^{\star}_{t})^\sfT (\bx^{\star}_{t+1} - \bx^{\star}_{t}) \nonumber \\
& \leq | - \bx_{t+1}^\sfT(\bx^{\star}_{t+1} - \bx^{\star}_{t}) + (\bx^{\star}_{t})^\sfT (\bx^{\star}_{t+1} - \bx^{\star}_{t}) | \\
& \leq \|  \bx_{t+1} \| \|\bx^{\star}_{t} - \bx^{\star}_{t+1}\| + \| \bx^{\star}_{t}\| \|\bx^{\star}_{t+1} - \bx^{\star}_{t}\| \\
& = \| \hat{\sfF}_t(\bx_{t}) \| \|\bx^{\star}_{t} - \bx^{\star}_{t+1}\| + \| \hat{\sfF}_t(\bx^{\star}_{t})\| \|\bx^{\star}_{t+1} - \bx^{\star}_{t}\| \\
& \leq 2 M_t \sigma_t  \label{eq:proof_afp_10}
\end{align}
\end{subequations}
Therefore, using~\eqref{eq:proof_afp_7},~\eqref{eq:proof_afp_9} in \eqref{eq:proof_afp_4} and~\eqref{eq:proof_afp_10}, one obtains the following bound:
\begin{align}
& \|\bx_{t+1} - \bx^{\star}_{t+1}\|^2   \leq \|\bx_{t} - \bx^{\star}_{t}\|^2 - \alpha_t(1 - \alpha_t) \|\sfT_t(\bx_{t}) - \bx_{t}\|^2  \nonumber \\
& + \alpha_t^2 e_{\sfT,t}^2  + 4 \alpha_t  e_{\sfT,t} M_t + 4 M_t \sigma_t + \sigma_t^2 \label{eq:proof_afp_11}
\end{align}
or, equivalently, 
\begin{align}
& \alpha_t(1 - \alpha_t) \|\sfT_t(\bx_{t}) - \bx_{t}\|^2 \nonumber \\
& \hspace{.8cm} \leq \|\bx_{t} - \bx^{\star}_{t}\|^2 - \|\bx_{t+1} - \bx^{\star}_{t+1}\|^2 \nonumber \\
& \hspace{.8cm} + \alpha_t e_{\sfT,t} ( 4M_t + \alpha_t e_{\sfT,t}) + \sigma_t (4 M_t + \sigma_t) \, .  \label{eq:proof_afp_12}
\end{align}
Summing~\eqref{eq:proof_afp_12} over $t = 1, 2, \ldots, T$ yields~\eqref{eq:fp_residual}.

\subsection{Proof of Corollary~\ref{cor:vanishing_error}}
\label{sec:proof_cor_vanishing_error}

Note that~\eqref{eq:sum_errors} implies that $\lim_{T \rightarrow \infty} \frac{1}{T}\sum_{t = 1}^{T} e_{\sfT,t}^2 = 0$ as:
\begin{align}
\sum_{t = 1}^T (e_{\sfT,t}/e_{\sfT})^2 \leq \sum_{t = 1}^T e_{\sfT,t}/e_{\sfT} \leq  o(T)/{e_{\sfT}}
\end{align}
implying that $\sum_{t = 1}^T e_{\sfT,t}^2 \leq e_{\sfT} o(T)$. Then,~\eqref{eq:fp_residual_3} can be shown from Theorem~\ref{thm:residual}. 

\subsection{Proof of Theorem~\ref{thm:convergence}}
\label{sec:proof_thm_convergence}

Bound $\|\bx_{t+1} - \bx^{\star}_{t+1}\|$ as: 
\begin{subequations}
\begin{align}
& \|\bx_{t+1} - \bx^{\star}_{t+1}\|  = \|\bx_{t+1} - \bx^{\star}_{t} - (\bx^{\star}_{t+1} - \bx^{\star}_{t})\| \label{eq_proof_ineq_1} \\
& \hspace{.2cm} = \|\hat{\sfF}_t(\bx_{t}) - \sfF_t(\bx^{\star}_{t}) - (\bx^{\star}_{t+1} - \bx^{\star}_{t})\| \label{eq_proof_ineq_2}  \\
& \hspace{.2cm} \leq \|\hat{\sfF}_t(\bx_{t}) - \sfF_t(\bx^{\star}_{t})\| + \sigma_t \label{eq_proof_ineq_3}  \\
& \hspace{.2cm} \leq \|{\sfF}_t(\bx_{t}) - {\sfF}_t(\bx^{\star}_{t})\| + \|\hat{\sfF}_t(\bx_{t}) - \sfF_t(\bx_{t})\| + \sigma_t \label{eq_proof_ineq_4}  \\
& \hspace{.2cm} \leq L_t \|\bx_{t} - \bx^{\star}_{t}\| + \alpha_t  \|\hat{\sfT}_t(\bx_{t}) - \sfT_t(\bx_{t})\| + \sigma_t \label{eq_proof_ineq_5}  \\
& \hspace{.2cm} \leq {L}_t \|\bx_{t} - \bx^{\star}_{t}\| + \alpha_t e_{\sfT,t} + \sigma_t \label{eq_proof_ineq_6}
\end{align}
\end{subequations}
where the definition of $\sigma_t$ was used in~\eqref{eq_proof_ineq_3}  and Assumption~\ref{asm:assump_L} was utilized to obtain~\eqref{eq_proof_ineq_5}. 
Therefore, 
\begin{align}
\|\bx_{t+1} - \bx^{\star}_{t+1}\| \leq L_t \|\bx_{t} - \bx^{\star}_{t}\| + \alpha_t e_{\sfT,t} + \sigma_t \, . \label{eq_proof_ineq_7}
\end{align}
Applying~\eqref{eq_proof_ineq_7} recursively for $\tau = 1, \ldots, t$ yields~\eqref{eq:bound_iter}.   

Next, take $\gamma := \sup_t \{\alpha_t e_{\sfT,t}\} + \sup_t \{\sigma_t\}$ and $L := \sup_t \{L_t\}$, where $L < 1$. Then,~\eqref{eq:bound_iter} is upper bounded by 
\begin{align}
\|\bx_{t+1} - \bx^{\star}_{t+1}\| & \leq \bar{c}^{(t,0)} \|\bx_1 - \bx^{\star}_1\| + \gamma \sum_{\tau = 1}^t \bar{c}^{(t, \tau)}   \label{eq:bound_iter_2} 
\end{align}
where $\bar{c}^{(t, \tau)} = 1$ is $\tau = t$ and $\bar{c}^{(t, \tau)} = L^{t-\tau+1}$ is $\tau = 1, \ldots, t-1$. The first term on the  right-hand-side of~\eqref{eq:bound_iter_2} vanishes with the increasing of $t$. The second term on the  right-hand-side is the sum of the first $t$ terms of a geometric series. Taking the limit for $t \rightarrow + \infty$ the result~\eqref{eq:fp_residual} follows.  

\subsection{Proof of Proposition~\ref{prop:gradient}}
\label{sec:proof_prop_gradient}

First, for each time $t$, $\nu \in (0, 2/K_t)$ then the fact that $\alpha_t = 1/(2 - \nu K_t/2)$ is proved in~\cite[Proposition 2.4]{Combettes15}. The exact and approximate maps $\sfT_t$ and $\hat{\sfT}_t$ can be expressed as: 
\begin{align}
\sfT_t(\bx) &= \frac{\alpha_t - 1}{\alpha_t} \bx + \frac{1}{\alpha_t}\textrm{proj}_{\cX_t}\{\bx - \nu \nabla f_t(\bx)\}
  \label{eq:proof_projgrad_1} \\
\hat{\sfT}_t(\bx) &= - \frac{1 - \alpha_t}{\alpha_t} \bx + \frac{1}{\alpha_t}\textrm{proj}_{\cX_t}\{\bx - \nu \by_t\}
  \label{eq:proof_projgrad_2} \, .
\end{align}
Therefore, using the non-expansive property of the projection operator, one has that:
\begin{align}
\|\sfT_t(\bx) - \hat{\sfT}_t(\bx) \| \leq \frac{\nu}{\alpha_t} \| \nabla f_t(\bx) - \by_t \|
  \label{eq:proof_projgrad_3} \, .
\end{align}
Using $\alpha_t = 1/(2 - \nu K_t/2)$ and the bound for $ \|\nabla f_t(\bx) - \by_t \|$, the result (i) follows. The result for (ii) builds on the strong convexity and strong smoothness of $f_t$; when $\nu \in (0, 2/K)$, then the operator $\sfI - \nu \nabla f_t$ is contractive, and the composition of a contractive operator and a non-expansive one is contractive~\cite{Combettes15}.

\bibliographystyle{IEEEtran}
\bibliography{biblio.bib,PaperCollectionAS.bib}

\vspace{.2cm}

\small{This work was authored in part by the National Renewable Energy Laboratory, operated by Alliance for Sustainable Energy, LLC, for the U.S. Department of Energy (DOE) under Contract No. DE-AC36-08GO28308.
The views expressed in the article do not necessarily represent the views of the DOE or the U.S. Government. The publisher, by accepting the article for publication, acknowledges that the U.S. Government retains a nonexclusive, paid-up, irrevocable, worldwide license to publish or reproduce the published form of this work, or allow others to do so, for U.S. Government purposes. }

 \end{document}